\documentclass[11pt,a4paper]{article}
\usepackage[latin1]{inputenc}
\usepackage[american]{babel}      % language
\usepackage[T1]{fontenc}
\usepackage{amsmath,amssymb,amsthm}
\usepackage{graphicx}
\usepackage{xcolor}   % pacchetto per scrivere a colori
\usepackage{geometry}
\usepackage{cancel}
\title{{\bf Simultaneous approximation by neural network operators with applications to Voronovskaja formulas}}
         
\author{ {\bf Marco Cantarini} and {\bf Danilo Costarelli}\thanks{Corresponding author} \\  
Department of Mathematics and Computer Science \\
            University of Perugia\\
        1, Via Vanvitelli, 06123 Perugia, Italy    \\  
{\small {\tt marco.cantarini@unipg.it} - {\tt danilo.costarelli@unipg.it}} }

\date{}

\newcommand{\mau}{\geq}
\newcommand{\miu}{\leq}

\newcommand{\N}{\mathbb{N}}
\newcommand{\R}{\mathbb{R}}
\newcommand{\Z}{\mathbb{Z}}

\newcommand{\be}{\begin{equation}}
\newcommand{\ee}{\end{equation}}
\newcommand{\phis}{\phi_{\sigma}}

\newcommand{\phil}{\phi_{\sigma_{\ell}}}

\newtheorem{definition}{Definition}[section]
\newtheorem{remark}[definition]{Remark}
\newtheorem{theorem}[definition]{Theorem}
\newtheorem{lemma}[definition]{Lemma}

\begin{document}

\maketitle 

\begin{abstract}
In this paper, we considered the problem of the simultaneous approximation of a function and its derivatives by means of the well-known neural network (NN) operators activated by sigmoidal function. Other than a uniform convergence theorem for the derivatives of NN operators, we also provide a quantitative estimate for the order of approximation based on the modulus of continuity of the approximated derivative.
%, including an asymptotic term depending on the decay behavior of the sigmoidal function and its derivatives.
Furthermore, a qualitative and quantitative Voronovskaja-type formula is established, which provides information about the high order of approximation that can be achieved by NN operators. To prove the above theorems, several auxiliary results involving sigmoidal functions have been established. At the end of the paper, noteworthy examples have been discussed in detail.
\end{abstract}

\section{Introduction} \label{sec1}

  The study of artificial neural networks (NNs) currently represents a very active field of research (see, e.g., \cite{GOL1,GUR1,BOGR1,EL1,DEHA1,Shen1,GRVO1}). 

Beyond their applications in various applied sciences, NNs are also used in Approximation Theory, as typical training processes with appropriate training set can be formally regarded as an approximation task using a specific (neural) method (see, e.g., \cite{CY,MAK1}). In this context, a considered research problem is the one of the so-called {\em simultaneous approximation}, namely the possibility to reconstruct both a given function $f$ and its derivatives, using sample/data taken only from the original function $f$, without requiring additional explicit information on its derivatives. It is clear that this kind of results showed the effectiveness of the considered reconstruction method of approximation that, just using one family of data allows reconstructing more than one analog signal, modeling real world situations. Some very preliminary and general results on the latter topic have been already established in \cite{LI1,XC1,XC2,SY1}.

In this area we can also collocate the theory of neural network (NN) operators, which provide one of the possible mathematical formalizations of the previous tools, characterized by a deterministic approach.

The theory of NN operators arose from the paper of Cardaliaguet and Euvrard \cite{CAEU1}, and it has been recently studied by several authors under different points of views. For some literature on this respect, the readers can see \cite{CACH1,CACH2,TUDU,BAKU1,QY1,KA1,WYZ1,AN1,BA1,BAAG,KA2,WYG1}.

The above quoted papers are all substantially based on the version of the NN operators introduced in \cite{COSP1}, which generalize and extend the ones of Cardaliaguet and Euvrard.

More precisely, in \cite{COSP1} the authors introduced discrete NN operators (denoted by $F_n$, $n \in \N^+$) activated by sigmoidal functions (\cite{CY}), i.e., by functions having limit 1, as $x \to +\infty$, and limit 0, as $x \to -\infty$. The main peculiarity of these operators is that they are suitable in order to get constructive approximation results; indeed, their expression is known a priori: they have $n(b-a)+1$ neurons (when they are considered on the interval $[a,b]$, with $a$, $b$ integers), they have fixed positive integer weights, and their coefficients can be effectively computed in term of discrete values of $f$. 

  For these operators pointwise convergence was obtained at continuity points, while the uniform convergence was obtained in the approximation of continuous functions on bounded intervals. Subsequently, also the multivariate theory has been developed in \cite{COSP2}, as one is expecting in case of neural network type models. Recently, also a "deep" version of the NN operators has been introduced in \cite{CO27} and also studied in \cite{SS1}.

A classical problem in Approximation Theory is the ''simultaneous approximation issue'' described above. To address this, it is essential to understand when derivatives of a given function can be approximated by the derivatives of a given family of approximation operators. We stress again that one of the main advantages of this kind of results is that, one can approximate the $s$-th derivative of a given function $f$ knowing only a discrete family of values of $f$, without any additional information on the values assumed by $f^{(s)}$. Even if this research problem has been addressed in the past from a very general point of view in the theory of NNs, for the case of NN operators it was not never studied yet.

Hence, in the present paper, we face the above issue for the mentioned NN operators activated by sigmoidal functions; a crucial question in such an approximation problem is provided by the computation of the asymptotic behavior of the so-called {\em truncated algebraic moments} of both the density functions $\phis$ (generated by suitable sigmoidal function $\sigma$) and its derivatives $\phis^{(s)}$. In order to do this, we have to require certain conditions on $\phis$, that are known in the literature as the Strang-Fix type condition (\cite{STFI1}). 

  The necessity to assume such a condition is not surprising, since this is very typical in Approximation Theory and Fourier Analysis when one deals with the problem of simultaneous approximation and/or high order of convergence by Voronovskaja type formulas.

  The proposed proof of the above result of simultaneous approximation is constructive, and also allows us to deduce a quantitative estimate for the order of approximation thanks to the use of the well-known modulus of continuity (see, e.g., \cite{DELO1}).

Finally, using the previous results concerning the truncated algebraic moments of $\phis$ we can also derive the so-called Voronovskaja type formulas, which provide information about the exact high order of approximation that can be achieved by NN operators. 

Finally, at the end of the paper a detailed discussion concerning examples of sigmoidal functions for which the above results hold have been provided. In particular, the case of the logistic function has been analyzed, and other important examples have been mentioned, such as  the case of sigmoidal functions generated by B-splines, from which we can deduce results for the celebrated ReLU activation function (see \cite{GO1,CYL1}).

In the latter examples, in order to prove that the previously mentioned Strang-Fix type assumption is satisfied, we have to resort to the application of the celebrated Poisson summation formula, involving the Fourier transform of the involved density functions, and their derivatives.

%%%%%%%

\section{Main definitions} \label{sec2}

In the paper, by the symbol $C(\Omega)$, with $\Omega=[a,b]$ or $\Omega=\R$, we denote the spaces of all uniformly continuous and bounded functions. Moreover, $C^m(\Omega)$, $m \in \N^+$, is the subspace of $C(\Omega)$ of functions $f$ having derivatives $f^{(s)}$, $s=1,...,m$, each one belonging to $C(\Omega)$. 

In what follows, to denote the derivative of a function or of an operator we also use the standard symbol $d^s/dx^s$.

A measurable function $\sigma: \R \to \R$ is called a {\em sigmoidal function} if $\lim_{x \to -\infty}\sigma(x)=0$ and $\lim_{x \to +\infty}\sigma(x)=1$.

   From now on, we always consider non-decreasing sigmoidal functions $\sigma$, satisfying the following assumptions:
\begin{itemize}
\item[$(\Sigma 1)$] $\sigma(x)-1/2$ is an odd function;
\item[$(\Sigma 2)$] $\sigma \in C^2(\R)$ is concave for $x \mau 0$;
\item[$(\Sigma 3)$] $\sigma(x)={\cal O}(|x|^{-\alpha-1})$ as $x \to -\infty$, for some $\alpha>0$.
\end{itemize}

  Now, we recall the definition of the density (kernel) function $\phis$ generated by $\sigma$:
$$
\phi_{\sigma}(x)\, :=\, \frac{1}{2}[\sigma(x+1)-\sigma(x-1)], \hskip1cm x \in \R.
$$

Clearly, based on the above assumptions, it is immediate to observe that $\phis$ is an even function taking
non-negative values, non-decreasing on $(-\infty ,0]$ and non-increasing on $%
[0,\infty )$. In addition to that, $\phis$ satisfies a list of
properties (see \cite{COSP1}) from which it
is enough here to recall that:
\be \label{series-equal-one}
\sum\limits_{k\in \mathbb{Z}}\phi _{\sigma }(x-k) = 1, \quad x\in 
\R;
\ee
and
\be 
1\ \geq\ \sum\limits_{k=\left\lceil na\right\rceil }^{\left\lfloor
nb\right\rfloor }\phi _{\sigma }(nx-k)\ \geq\ \phi _{\sigma }(1)\ >\ 0,
\ee
for every $a,b \in \R$ with $a<b$, $x\in \lbrack a,b]$ and $n\in \mathbb{N}_{+}$
such that $\left\lceil na\right\rceil \leq \left\lfloor nb\right\rfloor $,
where $\left\lceil t\right\rceil$ and $\left\lfloor
t\right\rfloor$ are the ceiling and the integer part of $t$, respectively, with $t\in \mathbb{R}$; and
\be
\phi _{\sigma}(x)=\mathcal{O}\left( \left\vert x\right\vert ^{-\alpha-1
}\right), \quad as \quad x \to +\infty,
\ee
where $\alpha >0$ is the parameter of $(\Sigma 3)$, from which we get $\phi_{\sigma} \in L^1({\mathbb{R}})$.

From the sake of completeness, we highlight that the term "kernel" arise from the classical theories of singular integrals and cardinal series, when we deal with functions satisfying the continuous or discrete versions of the usual partition of the unity property (\ref{series-equal-one}).

Let now $f:[a,b]\rightarrow \mathbb{R}$ be a bounded function, and $n\in 
\mathbb{N}_{+}$, such that $\left\lceil na\right\rceil \leq \left\lfloor
nb\right\rfloor$. The neural network (NN) operators activated by the
sigmoidal function $\sigma $, are defined as (see Definition 3.1 in \cite{COSP1})
\begin{equation} \label{NN}
F_{n}(f,x)=\frac{\displaystyle \sum\limits_{k=\left\lceil na\right\rceil }^{\left\lfloor
nb\right\rfloor } f\left( \frac{k}{n}
\right) \phi _{\sigma }\left( nx-k\right)  }{\displaystyle \sum\limits_{k=\left\lceil na\right\rceil }^{\left\lfloor
nb\right\rfloor }\phi _{\sigma }\left( nx-k\right) }\text{,} \quad x \in
I:=[a,b].  
\end{equation}%
Note that, it is known that the above operators are well-defined, e.g., in case of bounded functions (\cite{COSP1}).

  We now recall a pointwise and uniform convergence theorem for the above family of operators, which turn out to be, obviously, well-defined. 
\begin{theorem}[\cite{COSP1}]  
  Let $f: I \to \R$ be bounded. Then,
$$
     \lim_{n \to +\infty} F_n(f,x) = f(x),
$$
at each point $x \in I$ where $f$ is continuous. Moreover, if $f$ is continuous on the whole $I$, we have:
$$
 \lim_{n \to +\infty} \|F_n(f, \cdot) - f(\cdot) \|_{\infty} = 0.
$$
\end{theorem}

%%%

\section{Auxiliary results} \label{sec3}

In order to establish results concerning the simultaneous approximation of a function and its derivatives, we have to introduce the following notation.

 For any function $\Phi: \R \to \R$, and $\nu \in \N$, we can define the truncated algebraic moment of order $\nu$ (associated to the interval $[a,b]$, see \cite{COVI2}) by:
$$
m^n_{\nu}(\Phi, u)\ :=\ \sum_{k=\lceil na \rceil}^{\lfloor nb \rfloor}\Phi(u-k)\, (k-u)^{\nu}, \hskip1cm u \in \R.
$$ 
Moreover, we can define also the algebraic moment of order $\nu \in \N$ of $\Phi$, as:
$$
{\cal A}_\nu(\Phi, u)\ :=\ \sum_{k \in \Z} \Phi(u-k)\, (k-u)^{\nu}, \hskip1cm u \in \R,
$$
and the discrete absolute moments of $\Phi$, by:
$$
M_{\nu}(\Phi)\ :=\ \sup_{u \in \R}\, \sum_{k \in \Z}|\Phi(u-k)|\, |u-k|^{\nu}, \quad \quad \nu\mau 0.
$$ 
Note that, the terms "truncated algebraic moment" arise from the classical notions used in the fundamental book of Approximation Theory written by P.L. Butzer and R. Nessel (\cite{BUNE}); it is clear that by the symbol $m^n_{\nu}(\Phi, u)$ we denoted a truncated version of the ${\cal A}_\nu(\Phi, u)$, that are widely used in the theory of sampling and related issues.

We can state the following lemma which provides some useful properties of the discrete absolute moments of the density functions $\phis$.
\begin{lemma}[\cite{COVI2}] \label{lemma2}
Let $\sigma(x)$ be a fixed sigmoidal function. Then the following assertion holds:
$M_{\nu}(\phis)<+\infty$, for every $0 \miu \nu < \alpha$.
\end{lemma}

  Now, we have to prove an auxiliary result concerning the truncated algebraic moments of $\phis$. In order to do this, we assume the following additional condition on $\phis$:
\begin{itemize}

\item[$(\Sigma 4)$] $\sigma \in C^m(\R)$, $m \mau 2$, such that, there exist (sufficiently large) positive constants $K_s$ and $C_s$, with:
\be \label{m-regular}
|\sigma^{(s)}(x)|\ \miu C_s\, |x|^{-\beta-1}, \quad  |x|  \mau K_s, \quad s=1, 2, ..., m,
\ee
where $\beta>m + 1$. 
\end{itemize}
\vskip0.2cm

It is clear that, under condition $(\Sigma 4)$, since $\phis^{(s)}(x)=\frac12[\sigma^{(s)}(x+1)-\sigma^{(s)}(x-1)]$, $x \in \R$, the same inequality given in (\ref{m-regular}) for $\sigma^{(s)}$ can be easily obtained also for $\phis^{(s)}$, that is:
\be
|\phis^{(s)}(x)|\ \miu C_s\, |x|^{-\beta-1}, \quad  |x|  \mau K_s, \quad s=1, 2, ..., m,
\ee
where, for simplicity and without losing of generality, we can use the same notations for the involved constants $K_s$ and $C_s$.
\vskip0.2cm

Note that, for sigmoidal functions satisfying $(\Sigma 4)$, and proceeding as in the proof of Lemma \ref{lemma2}, we have that:
\be \label{mom-der-sigma}
M_\nu(\sigma^{(s)})\ <\ +\infty, \quad 0\, \miu\, \nu\, <\, \beta, \quad s=1, ..., m.
\ee
Obviously, reasoning as above, we also have:
\be \label{mom-der-phis}
M_\nu(\phis^{(s)})\ <\ +\infty, \quad 0\, \miu\, \nu\, <\, \beta, \quad s=1, ..., m.
\ee

\vskip0.2cm

In order to prove the following result, from now on, we always denote by the symbol $I_\delta$ the following interval:
\be
I_\delta\ :=\ [a+\delta, b-\delta], 
\ee
for $\delta>0$ sufficiently small, such that the above set is non-trivial. Moreover, for the sake of simplicity, we also assume that both $a$ and $b$ are integers.\footnote{Obviously, what follows can also be proved in case of non-integer values of $a$ and $b$, introducing suitable and tedious restrictions of the parameter $\delta$.} Finally, from now on, we always consider sigmoidal functions $\sigma$ satisfying $(\Sigma i)$, $i=1,...,4$.
\vskip0.2cm

  Now, we can prove the following.
\begin{lemma} \label{first-lemma}
Let $\sigma$ be a given sigmoidal function, and let $\delta>0$ be a fixed sufficiently small parameter. Then:
$$
\lim_{n \to +\infty} m^n_0\left(\phis^{(s)},nx\right)\ :=\
\begin{cases}
1, \quad s=0,\\
0, \quad s=1, ..., m,
\end{cases}
$$
where the above limit holds uniformly with respect to $x \in I_\delta$. More precisely, we have:
$$
m^n_0\left(\phis^{(s)},nx\right)\ =\ {\cal O}(n^{-\beta-1}), \quad as \quad n \to +\infty, \quad \quad s=1,...,m,
$$
uniformly with respect to $x \in I_\delta$, where $\beta$ is the parameter of condition $(\Sigma 4)$.
\end{lemma}
\begin{proof}
Let $s \in \{0, 1, ..., m\}$ be fixed. Observing that $\phis^{(s)}(x)=\frac12 [\sigma^{(s)}(x+1)-\sigma^{(s)}(x-1)]$, we can write what follows:
$$
m^n_0(\phis^{(s)},nx)\ =\ \sum_{k=na}^{nb}\phis^{(s)}(nx-k)\, =\ \frac12\, \sum_{k=na}^{nb} [\sigma^{(s)}(nx - k+1)-\sigma^{(s)}(nx-k -1)]
$$
\be \label{useful-equality}
=\ \frac12\, \left\{  \sigma^{(s)}(n[x-a] +1)+ \sigma^{(s)}(n[x-a]) \right\}\, +\, \frac12\, \left\{  \sigma^{(s)}(n[x-b] +1)+ \sigma^{(s)}(n[x-b]) \right\}
\ee
because the sum is telescoping. Now, since $x \in I_\delta$, it is immediate to see that:
$$
x-a\, \mau\, \delta, \quad \mbox{and} \quad x-b\, \miu\, -\delta,
$$
from which we also have:
\be \label{estimate-for-uniformity}
n[x-a]\, \mau\, \delta\, n\, \mau\, K_s, \quad \mbox{and} \quad n[x-b]\, + 1 \miu\, -\delta\, n\, +\, 1 \miu - K_s,
\ee
for $n \in \N^+$ sufficiently large, where $K_s>0$ is the parameter arising from $(\Sigma 4)$. Thus, by the above considerations, using (\ref{useful-equality}) with $s=0$, and recalling that $\lim_{u\to+\infty}\sigma(u)=1$ and $\lim_{u\to-\infty}\sigma(u)=0$, we immediately have:
\be \label{first-part-case-a}
\lim_{n \to +\infty} m^n_0(\phis,nx)\ =\ 1.
\ee 
Note that, the above limit holds uniformly with respect to $x \in I_\delta$ in view of the inequalities in (\ref{estimate-for-uniformity}).

  Moreover, using (\ref{useful-equality}) with $s\mau1$ and by the property (\ref{m-regular}), we obtain:
\be \label{decay-beta}
|m^n_0(\phis^{(s)},nx)|\ \miu\ \frac12\, C_s\, \left[  {1 \over (n \delta + 1)^{\beta + 1}}\, +\, {2 \over (n \delta )^{\beta + 1}}\, +\, {1 \over (n \delta - 1)^{\beta + 1}}\right],
\ee
from which we get:
\be \label{proof-case1}
\lim_{n \to +\infty} m^n_0(\phis^{(s)},nx)\ =\ 0, \quad \quad s\, =\, 1,...,m,
\ee
uniformly with respect to $x \in I_\delta$. Note that, from (\ref{decay-beta}) immediately follows the second part of the thesis. This completes the proof.
\end{proof}

In order to study the limit behavior of the truncated algebraic moments of order $\nu>0$ of $\phis$, we have to assume the following additional assumption on $\sigma$. We require that $\sigma$ is such that:
\begin{enumerate}

\item[$(\Sigma 5)$] ${\cal A}_j(\phis,x) = A_{0, j} \in \R$, for every $x \in \R$, $j=1, ..., m$, where $m \in \N$ is the parameter of condition $(\Sigma 4)$.

\end{enumerate} 

The condition $(\Sigma 5)$ is a kind of Strang-Fix type assumption (see \cite{STFI1}), very typical in Approximation Theory and Fourier Analysis when one deal with the problem of simultaneous approximation and high order of convergence by Voronovskaja type formulas (see, e.g., \cite{CANCV4} and the references therein).
\vskip0.2cm

Now, we can prove the following.
\begin{lemma} \label{lemma-cal-A}
Let $\sigma$ be a sigmoidal function assumed as above and satisfying $(\Sigma 5)$. Then:
$$
{\cal A}_j(\phis^{(s)},x)\ =:\ A_{s, j}\ =\ 0, \quad x \in \R,
$$
for every $j=0, ..., s-1$, and
$$
{\cal A}_s(\phis^{(s)},x)\ =:\ A_{s, s}\ =\ s!, \quad x \in \R,
$$
where $s=1, ..., m$, being $m$ the parameter of condition $(\Sigma 4)$.
\end{lemma}
\begin{proof}
Let $s \in \{1,...,m\}$ be fixed and $j \in \{0, ..., s-1\}$. From conditions (\ref{series-equal-one}), $(\Sigma 4)$ and $(\Sigma 5)$ we can deduce that:
$$
{d^s \over dx^s} \left(  \sum_{k \in \Z} \phis(x-k)\, \left(k-x \right)^j  \right)\ =\ 0, \quad x \in \R.
$$
Differentiating the above series term-by-term (this is possible since all the series of the derivatives are absolutely and uniformly convergent on $I_\delta$ by the choose of $\beta$), and by using recursively the formula for the derivative of the product of two functions, we obtain:
$$
\hskip-6cm 0\ =\ {d^s \over dx^s}\left(  \sum_{k \in \Z} \phis(x-k)\, \left(k - x\right)^j  \right)
$$
\be \label{expr-useful-deriv}
=\ \sum_{\ell = 0}^j \binom{j}{\ell}  (s)_{j-\ell}\, (-1)^{j+\ell}\, \sum_{k \in \Z} \phis^{(\ell+s-j)}(x-k)\, \left(k-x \right)^{\ell},
\ee
where by the symbol $(s)_{j-\ell}$ we denote the falling factorial of $s$, that in general is defined, for every non-negative integer $r$, as follows:
\be
(s)_{r}\ :=\ \begin{cases}
s \cdot (s-1) \cdot ... \cdot (s-r+1), \quad r \mau 1,\\
1, \quad \quad r=0.
\end{cases}
\ee 
Hence, from (\ref{expr-useful-deriv}) we can deduce that:
\be \label{derivative-moments-etc}
\sum_{\ell = 0}^j \binom{j}{\ell}  (s)_{j-\ell}\, (-1)^{j+\ell}\, {\cal A}_{\ell}(\phis^{(\ell+s-j)},x)\ =\ 0.
\ee
Now, from (\ref{derivative-moments-etc}), if we set $j=0$ we immediately get:
$$
{\cal A}_{0}(\phis^{(s)},x)\ =\ 0, \quad \mbox{with} \quad s \mau 1,
$$
if $j=1$, and using the above equality, we get:
$$
{\cal A}_{1}(\phis^{(s)},x)\ =\ 0, \quad \mbox{with} \quad s \mau 2,
$$
and so on, with $j$ running between $0$ and $s-1$. 
Similarly, even from (\ref{derivative-moments-etc}) and by the relation $(s)_{s-\ell}\, \ell ! = s!$, for any $s=j$ we can write:
\begin{equation}
{\cal A}_{s}(\phis^{(s)},x)\ =\ s!\, (-1)^{s+1}\sum_{\ell = 0}^{s-1}\, \binom{s}{\ell}  {1 \over \ell!}\, (-1)^{\ell}\, {\cal A}_{\ell}(\phis^{(\ell)},x). \label{evaluation A_s}
\end{equation}
Now we prove by induction on $s$ that $\mathcal{A}_{s}\left(\phi_{\sigma}^{(s)},x\right)=s!$ for every $x\in\mathbb{R}$. We already noticed that $\mathcal{A}_{0}\left(\phi_{\sigma},x\right)=\sum_{k\in\mathbb{Z}}\phi_{\sigma}\left(x-k\right)=1$, so now we fix $s\geq1$. Then we get, by formula (\ref{evaluation A_s}) and the inductive assumption that:
$$
{\cal A}_{s}(\phi_{\sigma}^{(s)},x)\ =\ s!\,(-1)^{s+1}\sum_{\ell=0}^{s-1}\,\binom{s}{\ell}\frac{1}{\ell!}\,(-1)^{\ell}\,{\cal A}_{\ell}(\phi_{\sigma}^{(\ell)},x)=\ s!\,(-1)^{s+1}\sum_{\ell=0}^{s-1}\,\binom{s}{\ell}\,(-1)^{\ell}=s!,$$ 
by the binomial theorem. This completes the proof.
\end{proof}
Now, we are able to prove what follows concerning the truncated algebraic moments of $\phis$.
\begin{lemma} \label{second-lemma}
Let $\sigma$ be a sigmoidal function satisfying $(\Sigma i)$, $i=1,...,5$. In particular, we assume that $\alpha>m+1$, where $\alpha$ and $m$ are the parameters of conditions $(\Sigma 3)$ and $(\Sigma 4)$, respectively. Then, we have:
$$
m^n_j(\phis^{(s)},nx)\ =\ \begin{cases}
A_{0,j}\, +\, {\cal T}^1_{0,j,n}(x) + {\cal T}^2_{0,j,n}(x), \quad  s=0, \quad $j=1, ..., m$ \\
\\
{\cal T}^1_{s,j,n}(x) + {\cal T}^2_{s,j,n}(x), \quad \quad \quad \quad s=1, ..., m, \quad   0 \miu j < s,\\
\\
s!\, +\, {\cal T}^1_{s,j,n}(x) + {\cal T}^2_{s,j,n}(x), \quad \quad s=1, ..., m, \quad j=s,
\end{cases}
$$
where the above functions $\mathcal{T}^1_{s,j,n}$ and  $\mathcal{T}^2_{s,j,n}$ satisfy the following inequalities:
$$
|{\cal T}^1_{0,j,n}(x) + {\cal T}^2_{0,j,n}(x)|\ \miu\ 2\, \bar R_{\alpha, 0, j}\, n^{-(\alpha-j+1)/2},
$$
\be \label{resto-zero}
|{\cal T}^1_{s,j,n}(x) + {\cal T}^2_{s,j,n}(x)|\ \miu\ 2\, \bar R_{\beta, s, j}\, n^{-(\beta-j+1)/2},
\ee
for every $x \in I_\delta$, $n \in \N^+$ sufficiently large, for suitable absolute positive constants $\bar R_{\alpha, 0, j}$, $\bar R_{\beta, s, j}$ depending only on $\alpha$ and $\beta$, respectively ($\beta$ is the one of condition $(\Sigma 4)$), $s$ and $j$. 
\end{lemma}
\begin{proof}
Let $s \in \{ 0, 1, ..., m \}$ be fixed. Using Lemma \ref{lemma-cal-A} (where condition $(\Sigma 5)$ is one of the required assumption), for $j=0, 1, ..., s$, we can write what follows:
$$
\hskip-5cm m^n_j(\phis^{(s)},nx)\ =\ \sum_{k=na}^{nb}\phis^{(s)}(nx-k)\, (k-nx)^j\ 
$$
$$
=\ \sum_{k \in \Z} \phis^{(s)}(nx-k)\, (k-nx)^j\ - \sum_{k \miu na -1}\phis^{(s)}(nx-k)\, (k-nx)^j\ - \sum_{k \mau nb+1}\phis^{(s)}(nx-k)\, (k-nx)^j\
$$
$$
=\ A_{s, j} - \sum_{k \miu na -1}\phis^{(s)}(nx-k)\, (k-nx)^j\ - \sum_{k \mau nb+1}\phis^{(s)}(nx-k)\, (k-nx)^j\ =:\ A_{s, j} + {\cal T}^1_{s,j,n}(x) + {\cal T}^2_{s,j,n}(x).
$$
Now, recalling (\ref{estimate-for-uniformity}), for every $n > K_s / \delta$, where $K_s$ is the parameter arising from condition $(\Sigma 4)$, we get:
$$
\left| {\cal T}^1_{s,j,n}(x)  \right|\ \miu\ C_s\, \sum_{k \miu na -1}{ 1 \over |nx-k|^{\beta-j+1}}\ =\ C_s\, \sum_{i=1}^{+\infty}{ 1 \over \left\{ (n(x-a)+i)^2\right\}^{(\beta-j+1)/2}},
$$
(with $\alpha$ of condition $(\Sigma 3)$ in place of $\beta$ in the case $s=0$).
Since $x \in I_\delta$ we have:
$$
(n(x-a)+i)^2\ \mau\ (n \delta + i)^2\ =\ n^2\, \delta^2 + i^2 + 2 n \delta i\ \mau\ 2 n \delta i, \quad i \in \N^+,
$$
thus:
\be
\left| {\cal T}^1_{s,j,n}(x)  \right|\ \miu\ n^{-(\beta-j+1)/2}\frac{C_{s}}{(2\delta)^{(\beta-j+1)/2}}\,\zeta\left(\frac{\beta-j+1}{2}\right)\ =:\ n^{-(\beta-j+1)/2} \bar R_{\beta, s, j},
\ee
where $\zeta(s)$ is the Riemann-Zeta function, with $\bar R_{\beta, s, j} \in \R^+$ since $\beta>m+1$ and so $(\beta-j+1)/2 >1$, in view of condition $(\Sigma 4)$ (also $\bar R_{\beta, 0, j} \in \R^+$ since $\alpha>m+1$, hence $(\alpha-j+1)/2 >1$). Similarly, for $n$ sufficiently large, we have:
$$
\left| {\cal T}^2_{s,j,n}(x)  \right|\ \miu\ C_s \sum_{k \mau nb+1} {1 \over |nx-k|^{\beta-j+1}}\ =\ C_s\, \sum_{i=1}^{+\infty}{ 1 \over \left\{ (n(x-b)-i)^2\right\}^{(\beta-j+1)/2}},
$$
(again with $\alpha$ of condition $(\Sigma 3)$ in place of $\beta$ in the case $s=0$). Hence, for every $x \in I_\delta$:
$$
n(x-b)-i\ \miu\ -\delta n-i\ <\ 0, \quad i \in \N^+,
$$
from which we have:
$$
(n(x-b)-i)^2\ \mau\ (\delta n + i)^2\ \mau\ 2 n \delta i.
$$
Thus, proceeding as above we immediately obtain:
$$
\left| {\cal T}^2_{s,j,n}(x)  \right|\ \miu\ n^{-(\beta-j+1)/2} \bar R_{\beta, s, j},
$$
(again $\alpha$ in place of $\beta$ when $s=0$).
Finally, using the previous estimates for ${\cal T}^1_{s,j,n}(x)$ and ${\cal T}^2_{s,j,n}(x)$, together with assumption $(\Sigma 5)$ and Lemma \ref{lemma-cal-A}, the proof immediately follows.
\end{proof}

%%%

\section{Simultaneous approximation with applications to Voronovskaja formulas} \label{sec4}

In order to face the problem of simultaneous approximation of a given $f:I\to \R$ and its derivatives, in what follows we have to consider approximations only on intervals $I_\delta$.

Based on the results established in the previous section, under condition $(\Sigma 4)$ on $\sigma$, it is clear that, to work in $I_\delta$ we can simplify the definition of the NN operators $F_n$ as follows:
\be
\widetilde{F}_n(f, x)\ :=\ \sum_{k=na}^{nb}f\left({k \over n}\right)\, \phis(nx-k), \quad \quad x \in I_\delta,
\ee
where $f:I \to \R$.
\vskip0.2cm 

  The above claim is due to Lemma \ref{first-lemma} which asserts that the denominator of $F_n$, i.e., $\sum_{k=na}^{nb}\phis(nx-k)$ provide uniform approximation of the unitary constant function ${\bf 1}: I \to \R$, ${\bf 1}(x)=1$, $x \in I$, on $I_\delta$.
\vskip0.2cm

Hence, the new family of NN operators $\widetilde{F}_n$ converges uniformly to $f$ on $I_\delta$, for any continuous function $f:I \to \R$.
\vskip0.2cm

Indeed, one can write:
$$
\left| \widetilde{F}_n(f,x) - f(x)\right| \miu\ \left| \widetilde{F}_n(f,x) - f(x) \sum_{k=na}^{nb}\phis(nx-k)\right|\ +\ \|f\|_\infty \left| \sum_{k=na}^{nb}\phis(nx-k) - 1\right|
$$
$$
\miu\ \sum_{k=na}^{nb}\left| f\left( {k \over n} \right)  - f(x)\right| \phis(nx-k)\ +\ \|f\|_\infty \left| m^n_0(\phis,nx) - 1\right|, \quad \quad x \in I_\delta.
$$
Then, by the above inequality, using Lemma \ref{first-lemma} with $s=0$ and repeating the proof of Theorem 3.2 of \cite{COSP1}, we immediately get the convergence.

In practice, the denominator of $F_n$ represents a sort of correction term, useful to get the uniform convergence to the whole interval $I$ instead of $I_\delta$ only.

  From now on in the paper we decided to consider only the modified operators $\widetilde{F}_n$. Obviously, the following results can be naturally extended also to the original version of the NN operators $F_n$ restricted to $I_\delta$.

Before stating our main result, we recall the following basic definition of Approximation Theory. For any continuous function $f:I \to \R$, we define its modulus of continuity as follows:
$$
\omega(f,h)\ :=\ \sup \left\{ |f(x)-f(t)|:\ t,x \in I,\ |t-x|\miu h\right\}, \quad h>0.
$$
It is well-known that, in general:
$$
\lim_{h \to 0^+} \omega(f,h)\ =\ 0.
$$

We can now prove the following theorem of simultaneous approximation, in which we also provide a quantitative estimate.
\begin{theorem} \label{th-simultaneous}
Let $\sigma$ be a sigmoidal function satisfying conditions $(\Sigma i)$, $i=1,...,5$. In particular, we assume that the condition $(\Sigma 4)$ is satisfied with $m \in \N^+$ and $\beta>2m$.

If $f \in C^m(I)$ be fixed. Then:
$$
\lim_{n \to +\infty} {d^s \over dx^s}\left( \widetilde{F}_n(f,x) \right)\ =\ f^{(s)}(x),  
$$
uniformly with respect to $x \in I_\delta$, $s=1, ...,m$. In particular, the following quantitative estimate holds:
$$
\hskip-8cm \left|  {d^s \over dx^s}\left( \widetilde{F}_n(f,x) \right) - f^{(s)}(x)\ \right|
$$
$$
\miu\ 2\, \sum_{i=0}^{s}{\|f^{(i)}\|_\infty \over i!}\, \bar R_{\beta, s, j}\, n^{s - i -(\beta-i+1)/2} 
+\ {\omega(f^{(s)}, 1/n) \over s!}\, \left[ M_s(\phis^{(s)})\ +\ M_{s+1}(\phis^{(s)})  \right]<+\infty,
$$
$x \in I_\delta$, $n \in \N^+$ sufficiently large, where $\bar R_{\beta, s, j}$ are the positive constants arising from Lemma \ref{second-lemma}.
\end{theorem}
\begin{proof}
Let $s \in \{ 1, ..., m\}$ and $x \in I_\delta$ be fixed. Using for $f$ the Taylor formula of order $s$ with Lagrange remainder, we have:
$$
f(k/n)\ =\ \sum_{i=0}^s {f^{(i)}(x) \over i!}\left({k \over n} - x\right)^i\ +\ {1 \over s!}\, \left[ f^{(s)}(\theta_{k/n,x}) - f^{(s)}(x) \right] \left({k \over n} - x\right)^s,
$$
where $\theta_{k/n,x}$ is a suitable number between $x$ and $k/n$, $k=na, ..., nb$. Replacing the above formula in the definition of ${d^s \over dx^s}\left( \widetilde{F}_n(f,x) \right)$ and using Lemma \ref{second-lemma} with $s=1,...,m$, we can write what follows:
$$
{d^s \over dx^s}\left( \widetilde{F}_n(f,x) \right)\ =\ n^s \sum_{k=na}^{nb}f\left({k \over n}\right)\, \phis^{(s)}(nx-k)\  
$$
$$
=\ n^s \sum_{k=na}^{nb}\left[  \sum_{i=0}^s {f^{(i)}(x) \over i!}\left({k \over n} - x\right)^i\ +\ {1 \over s!}\, \left[ f^{(s)}(\theta_{k/n,x}) - f^{(s)}(x) \right] \left({k \over n} - x\right)^s  \right]\, \phis^{(s)}(nx-k)
$$
$$
=\ n^s\, \sum_{i=0}^s{f^{(i)}(x) \over i!} \sum_{k=na}^{nb} \left({k \over n} - x\right)^i \phis^{(s)}(nx-k)
$$
$$
+\ n^s\, {1 \over s!}\,  \sum_{k=na}^{nb} \left[ f^{(s)}(\theta_{k/n,x}) - f^{(s)}(x) \right] \left({k \over n} - x\right)^s  \phis^{(s)}(nx-k)
$$
$$
=\ n^s\, \sum_{i=0}^s{f^{(i)}(x) \over i!}\, n^{-i}\,  m^n_i(\phis^{(s)},nx) +\ {n^s \over s!}\,  \sum_{k=na}^{nb} \left[ f^{(s)}(\theta_{k/n,x}) - f^{(s)}(x) \right] \left({k \over n} - x\right)^s  \phis^{(s)}(nx-k)
$$
$$
=\ n^s\, \sum_{i=0}^{s}{f^{(i)}(x) \over i!}\, n^{-i}\, \left[ {\cal T}^1_{s,i,n}(x) + {\cal T}^2_{s,i,n}(x) \right]\ +\ f^{(s)}(x) 
$$
$$
+\ {n^s \over s!}\,  \sum_{k=na}^{nb} \left[ f^{(s)}(\theta_{k/n,x}) - f^{(s)}(x) \right] \left({k \over n} - x\right)^s  \phis^{(s)}(nx-k).
$$
Concerning the last term of the above equality we can note that:
$$
{n^s \over s!}\,  \left| \sum_{k=na}^{nb} \left[ f^{(s)}(\theta_{k/n,x}) - f^{(s)}(x) \right] \left({k \over n} - x\right)^s  \phis^{(s)}(nx-k)   \right|
$$
$$
\miu\ {1 \over s!}\, \sum_{k=na}^{nb} \omega\left(f^{(s)}, |\theta_{k/n,x} - x|\right) |k-nx|^s\, |\phis^{(s)}(nx-k)|.
$$
Since:
$$
|\theta_{k/n,x} - x|\ \miu\ \left| {k \over n}-x  \right|, 
$$
and we know that the general inequality:
$$
\omega(f^{(s)}, \lambda h)\ \miu (1+\lambda)\, \omega(f^{(s)}, h), \quad \lambda,\, h >0,
$$
holds (see, e.g., \cite{DELO1}), we obtain:
$$
{n^s \over s!}\,  \left| \sum_{k=na}^{nb} \left[ f^{(s)}(\theta_{k/n,x}) - f^{(s)}(x) \right] \left({k \over n} - x\right)^s  \phis^{(s)}(nx-k)   \right|
$$
$$
\miu\ {\omega(f^{(s)}, 1/n) \over s!}\, \sum_{k=na}^{nb} \left( 1\, +\, |k-nx|  \right)\, |k-nx|^s\, |\phis^{(s)}(nx-k)|
$$
$$
\miu\ {\omega(f^{(s)}, 1/n) \over s!}\, \left[ M_s(\phis^{(s)})\ +\ M_{s+1}(\phis^{(s)})  \right]\ <\ +\infty.
$$
Note that, the constant $M_s(\phis^{(s)}) + M_{s+1}(\phis^{(s)})$ is finite in view of (\ref{mom-der-phis}) and since $(\Sigma 4)$ holds for $\beta>2m \mau m+1 \mau s+1$. 

In conclusion, using the same notations considered in the proof of Lemma \ref{second-lemma}, and what we have written above, we obtain:
$$
\hskip-6cm \left|  {d^s \over dx^s}\left( \widetilde{F}_n(f,x) \right) - f^{(s)}(x)\ \right|
$$
$$
\miu\ 2\, \sum_{i=0}^{s}{\|f^{(i)}\|_\infty \over i!}\, \bar R_{\beta, s, i}\, n^{s - i -(\beta-i+1)/2} 
+\ {\omega(f^{(s)}, 1/n) \over s!}\, \left[ M_s(\phis^{(s)})\ +\ M_{s+1}(\phis^{(s)})  \right]\ \longrightarrow\ 0,
$$
as $n \to +\infty$, where $s - i -(\beta-i+1)/2<0$, when $i=0,...,s$, since we assumed $\beta>2m \mau 2m-1$.
\end{proof}

By an approach similar to that one used above, also Voronovskaja-type formulas can be derived. In general, Voronovskaja formulas provide information concerning the exact order of approximation of a given family of approximation operators.  We can prove the following.
\begin{theorem} \label{th-voronovskaja}
Let $\sigma$ be a sigmoidal function satisfying conditions $(\Sigma i)$, $i=1,...,5$. In particular, we assume that the parameter $\alpha>2m$, where $\alpha$ and $m$ are the parameters of conditions $(\Sigma 3)$ and $(\Sigma 4)$, respectively. Moreover, we also assume that the condition $(\Sigma 5)$ is satisfied for $A_{0,j}=0$, $j=1,...,m-1$ and $A_{0,m}\neq 0$. If $f \in C^m(I)$ be fixed. Then:
$$
\lim_{n \to +\infty} n^{m}\left[ \widetilde{F}_n(f,x)  \, -\ f(x)\, \right]\ =\ {f^{(m)}(x)\, A_{0,m} \over m!}, \quad x \in I_\delta.  
$$
In particular, the following quantitative estimate holds:
$$
\left|n^{m}\left[ \widetilde{F}_n(f,x)  \, -\ f(x)\, \right]\ -\  {f^{(m)}(x) \over m!} A_{0,m} \right| \miu\ 2 \sum_{i=0}^m{\|f^{(i)}\|_\infty \over i!}\, \bar R_{\alpha, 0, j}\, n^{m-i - (\alpha - i +1)/2}
$$
$$
+\ {1 \over m!}\, \omega(f^{(m)},1/n)\, \left[ M_{m+1}(\phis) + M_m(\phis)  \right]\ <\ +\infty,
$$
$x \in I_\delta$, $n \in \N^+$ sufficiently large, where the constants $R_{\alpha, 0, j}$ are the ones arising from Lemma \ref{second-lemma}.
\end{theorem}
\begin{proof}
Proceeding as in the proof of Theorem \ref{th-simultaneous}, using the Taylor formula of order $m$ with Lagrange remainder for $x \in I_\delta$, we immediately have:
$$
\hskip-5cm \widetilde{F}_n(f,x)\ =\ 
\sum_{i=0}^m{f^{(i)}(x) \over i!}\, n^{-i}\,  m^n_i(\phis,nx)\
$$
$$
+\ {1 \over m!}\,  \sum_{k=na}^{nb} \left[ f^{(m)}(\theta_{k/n,x}) - f^{(m)}(x) \right] \left({k \over n} - x\right)^m  \phis(nx-k),
$$
where $\theta_{k/n,x}$ are suitable values between $x$ and $k/n$. Using Lemma \ref{second-lemma} in the case $s=0$, together with the fact that $A_{0,0}=1$ (by (\ref{series-equal-one})), $A_{0,j}=0$, for $j=1,...,m-1$, and $A_{0,m} \neq 0$, we obtain:
$$
\widetilde{F}_n(f,x)\ =\ f(x)\ +\ 
\sum_{i=0}^m{f^{(i)}(x) \over i!}\, n^{-i}\,  \left[ {\cal T}^1_{0,i,n}(x) + {\cal T}^2_{0,i,n}(x) \right]\ +\ n^{-m}\, {f^{(m)}(x) \over m!} A_{0,m}
$$ 
$$
+\ {1 \over m!}\,  \sum_{k=na}^{nb} \left[ f^{(m)}(\theta_{k/n,x}) - f^{(m)}(x) \right] \left({k \over n} - x\right)^m  \phis(nx-k).
$$
Rearranging the above terms we have:
$$
n^{m}\left[ \widetilde{F}_n(f,x)  \, -\ f(x)\, \right]\ -\  {f^{(m)}(x) \over m!} A_{0,m} =\ \sum_{i=0}^m{f^{(i)}(x) \over i!}\, n^{m-i}\,  \left[ {\cal T}^1_{0,i,n}(x) + {\cal T}^2_{0,i,n}(x) \right]
$$
$$
+\ {1 \over m!}\,  \sum_{k=na}^{nb} \left[ f^{(m)}(\theta_{k/n,x}) - f^{(m)}(x) \right] \left(k-nx\right)^m  \phis(nx-k).
$$
Now, recalling the estimates for the terms ${\cal T}^1_{0,j,n}(x)$ achieved in Lemma \ref{second-lemma}, and estimating the above remainder term as made in the proof of Theorem \ref{th-simultaneous}, we immediately obtain:
$$
\left|n^{m}\left[ \widetilde{F}_n(f,x)  \, -\ f(x)\, \right]\ -\  {f^{(m)}(x) \over m!} A_{0,m} \right| \miu\ 2 \sum_{i=0}^m{\|f^{(i)}\|_\infty \over i!}\, \bar R_{\alpha, 0, i}\, n^{m-i - (\alpha - i +1)/2}
$$
$$
+\ {1 \over m!}\, \omega\left(f^{(m)},1/n\right)\, \left[ M_{m+1}(\phis) + M_m(\phis)  \right]\ <\ +\infty.
$$
This completes the proof.
\end{proof}
\begin{remark} \rm
Note that, the first part of Theorem \ref{th-voronovskaja} replaces the results established in Theorem 3.1 and Theorem 3.3 of \cite{COVI2}. Indeed, the above Voronovskaja formula holds only on $I_\delta$, and it can not be valid in the whole $I$ as instead erroneously claimed in \cite{COVI2}. Furthermore, the quantitative estimate established in Theorem \ref{th-voronovskaja} is completely new.
\end{remark}

\begin{remark} \rm
Note that, it is well-known that in general a Voronovskaja formula gives an exact indication in what is the exact and also the best possible (saturation) order of approximation of a given family of approximation operators. To this purpose, we can mention the pioneer and classical work provided in 1912 by Voronovskaja in the case of Bernstein operators, and also the very general theorems established by Nishishirao on the saturation order of families of linear operators  (\cite{NI1}).  Here, Voronovskaja formulas play a very central (and crucial) role in determining the best possible order of approximation.
\end{remark}

\section{Particular cases and examples} \label{sec5}

In this section we discuss some noteworthy examples of sigmoidal functions for which we can apply the previous approximation results.
\vskip0.2cm

 We begin with the case of the NN operators activated by the logistic function (see e.g., \cite{COSA1,COCOKA1,CO27}), that is defined as follows:
$$
\sigma_{\ell}(x)\ :=\ \left( 1 + e^{-x} \right)^{-1}, \hskip1cm x \in \R.
$$
It is well-known that the logistic function (or sigmoidal functions, in general) have been considered in such models in view of their possible interpretation as activation function of the biological neurons, that is, the function that determines the activation or the non-activation of the artificial neuron. NN operators activated by logistic functions have been widely studied, see e.g. \cite{CSAMBV1,CO1}. 

Clearly, $\sigma_{\ell}$ is a smooth function and it satisfies all the assumptions $(\Sigma\, {\rm i})$, ${\rm i}=1,2,3$, of Section \ref{sec2}. Indeed, $(\Sigma 1)$ easily follows since $\sigma_{\ell}(x)-\frac12=\frac12-\sigma_{\ell}(-x)$, $x \in \R$, while $(\Sigma 2)$ follows computing the second derivative of $\sigma_\ell$ (see \cite{COSP1} for mode details). In particular, due to its exponential decay to zero as $x \to -\infty$, condition $(\Sigma 3)$ is satisfied for every $\alpha >0$, hence in view of Lemma \ref{lemma2}, it turns out that $M_{\nu}(\phi_{\sigma_{\ell}}) < +\infty$, for every $\nu \mau 0$. Similarly, also any derivative of $\sigma_{\ell}$ has exponential decay hence for $x \to \pm \infty$, hence also condition $(\Sigma 4)$ holds with any $m \in \N^+$. 

It is clear that, from what we proved, a very crucial issue is represented by the validity of condition $(\Sigma 5)$. In order to understand if such property is satisfied, the following classical result of Fourier Analysis is crucial. 

We refer to the celebrated Poisson summation formula (see, e.g., \cite{BUNE}) that claims that, for any continuous $\Phi: \R \to \R$, such that $g(u):=(-i u)^{\nu} \Phi(u)$, $u \in \R$, belongs to $L^1(\R)$, for some $\nu \in \N$ (where $i$ denotes the complex unit), the following equality holds: 
$$
(- i)^{\nu}\, \sum_{k \in \Z} (u-k)^{\nu} \Phi(u-k)\  =\ \sum_{k \in \Z}(\widehat{\Phi})^{(\nu)} (2 \pi k)\, e^{i 2 \pi k u}, \hskip0.8cm u \in \R,
$$
where $\widehat{\Phi}(v)=\int_{\R}\Phi(t)\, e^{-i v t}\, dt$, $v \in \R$, denotes the usual Fourier transform of $\Phi$, and $(\widehat{\Phi})^{(\nu)}$ its derivative of order $\nu$. 

In practice, the above formula allows to compute the algebraic moment ${\cal A}_\nu(\Phi,u)$, $u \in \R$, of a given density function, knowing the values assumed by its Fourier transform (and its derivatives) at the nodes $2 \pi k$, $k \in \Z$. 
\vskip0.2cm

\begin{figure}[h!]
\centering
\includegraphics[scale=0.5]{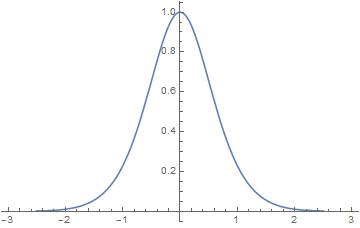}
\hskip0.5cm
\includegraphics[scale=0.5]{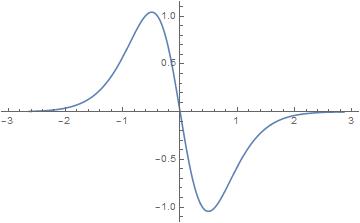}
\caption{The Fourier transform of $\phil$ (on the left) and its first derivative (on the right).} \label{fig1}
\end{figure}

Concerning the logistic function $\sigma_{\ell}(x)$, we can firstly observe that its corresponding density function $\phil(x)$ is band-limited (see, e.g., \cite{COVI2}), i.e., its Fourier transform has compact support. More precisely, we are able to show, using the CAS Mathematica, that $\widehat{\phil}$ has the following analytical expression:
$$
\widehat{\phil}(v)\ :=\ \frac{e^{-(i+\pi) v}}{2}\, \left\{ \beta\left(-1/e;\, 1-i v,\, 0\right)\, +\, e^{2(i+\pi)v}\, \beta\left(-1/e;\, 1+i v,\, 0\right)   \right.
$$
\be
\hskip-0.8cm -\ \left. e^{2 i v}\, \beta\left(-e;\, 1-i v,\, 0\right)\, -\, e^{2\pi v}\, \beta\left(-e;\, 1+i v,\, 0\right) \right\},
\ee
$v \in \R$, where here $\beta\left(x;\, y,\, z \right)$ denotes the usual incomplete Euler $\beta$-function, defined by:
$$
\beta\left(x;\, y,\, z \right)\ :=\ \int_0^x t^{y-1}\, (1-t)^{z-1}\, dt.
$$
Now, observing that $supp\, \widehat{\phil} \subset [-3, 3]$ (see Fig. \ref{fig1}), we can deduce that also $supp\, \widehat{\phil}^{(\nu)} \subset [-3, 3]$, $\nu \in \N^+$, and this implies that $\widehat{\phil}(2 \pi k) = \widehat{\phil}^{(\nu)}(2 \pi k) = 0$, for every $k \in \Z \setminus \left\{0 \right\}$, for every $\nu \in \N^+$. 

 Then, by applying the above Poisson summation formula, we can obtain that:
\be \label{momenti-logisitica}
{\cal A}_{\nu}(\phil, x)\ =\ i^{-\nu}(\widehat{\phil})^{(\nu)} (0)
\ee
$x \in \R$, $\nu \in \N$.  Since it is easy to verify that $\widehat{\phil}(0) = 1$, $(\widehat{\phil})' (0)=0$, $(\widehat{\phil})'' (0) \approx -3.6232$, and so on, we immediately have that all the algebraic moments ${\cal A}_{\nu}(\phil, x)$ are constants, $\nu \in \N$, hence Theorem \ref{th-simultaneous} holds, for every $m \in \N^+$. In practice, by the $s$-th derivative of $\widetilde{F}_n f$ activated by $\sigma_{\ell}$, we are able to approximate $f^{(s)}$, $s \mau 1$, and the order of approximation is determined by $\omega(f^{(s)}, 1/n)$. In particular, noting that:
\be \label{moments-logistic-voronov}
{\cal A}_{0}(\phil, x)\ =\ 1, \quad {\cal A}_{1}(\phil, x)\ =\ 0, \quad {\cal A}_{2}(\phil, x) \approx 3.6232 \neq 0, \quad x \in \R,
\ee
it is clear that in the case of the logistic function, the Voronovskaja formula of Theorem \ref{th-voronovskaja} holds with $m=2$. This result is perfectly coherent with what we know from the classical Korovkin theory for positive linear operators (see \cite{ALCA}), which asserts that the order $2$, is the best possible that can be achieved in these cases. 
\vskip0.2cm

The above computations and remarks can be given for other useful examples of sigmoidal functions. On this respect, we can apply the results established in this paper also for the hyperbolic tangent sigmoidal function (see \cite{COSP2,COVI21,COCO1}, and for its Fourier transform see \cite{COVI2}, Section 5.2), or by sigmoidal functions generated by the central B-spline of order $n$ (see \cite{CACOVI1,CACOCOGAVI1,CCNP1}, and for its Fourier transform see \cite{COVI2}, Section 5.3).

In the latter case, also some applications to the case of NN operators based on the ReLU functions can be easily deduced. For more details, see \cite{CO1,CO27,CP1}.

\begin{remark} \rm
Note that, the approximation results established in the present paper can suggest some possible applications in the field of signal/image reconstruction. Indeed, the theoretical results (both the simultaneous approximation and the Voronovskaja theorems) here established suggest that the considered approximation tools have a good accuracy in reconstructing data. Hence, the implementation of the multivariate version of such operators can be used to accurately reconstruct images (in fact providing an efficient algorithm for image resizing) and also to edge detection. In particular, the latter application is exactly the case in which the simultaneous approximation provides its better performance. Indeed, the (discrete) derivative of an image emphasizes the border of the figures of a given image; this means that by the derivatives of the NN operators one can detect all the border's figures. This can result in possible useful application in the setting of biomedical imaging, in which the detection of the border can produce automated diagnosis. In this sense we refer to the cancer detection from biomedical images, such as CT and/or RM images. 
\end{remark}

\section{Final remarks, conclusions and future developments}

In the present paper we study the problem of the simultaneous approximation and of the Voronoskaja formulas for a family of NN operators. In order to prove such theorems a very crucial point is the establishing of auxiliary results involving the density functions $\phis$ and their derivatives. Moreover, also than practical examples (see Section \ref{sec5}), some possible extensions for considering real world problems in biomedical diagnosis have been discussed. As future works, the results given in the present paper can be extended to the case of the recently introduced deep NN operators (see \cite{CO1}), or they can also be extended to the case of $L^p$-norm approximation, e.g., for approximating functions belonging to Sobolev spaces, or to establish estimates using different moduli of smoothness, such as the averaged moduli of smoothness.

%%%%

\section*{Acknowledgments}

{\small The authors are members of the Gruppo Nazionale per l'Analisi Matematica, la Probabilit\`a e le loro Applicazioni (GNAMPA) of the Istituto Nazionale di Alta Matematica (INdAM), of the network RITA (Research ITalian network on Approximation), and of the UMI (Unione Matematica Italiana) group T.A.A. (Teoria dell'Approssimazione e Applicazioni). 
}

\section*{Funding}

{\small The authors have been partially supported within the (1) 2024 GNAMPA-INdAM Project "Tecniche di approssimazione in spazi funzionali con applicazioni a problemi di diffusione" (CUP E53C23001670001) (2) 2024 GNAMPA-INdAM Project "Ricostruzione di segnali, tramite operatori e frame, in presenza di rumore" (CUP E53C23001670001), (3) "National Innovation Ecosystem grant ECS00000041 - VITALITY", funded by the European Union - NextGenerationEU under the Italian Ministry of University and Research (MUR), and (4) PRIN 2022 PNRR: "RETINA: REmote sensing daTa INversion with multivariate functional modeling for essential climAte variables characterization", funded by the European Union under the Italian National Recovery and Resilience Plan (NRRP) of NextGenerationEU, under the  MUR (Project Code: P20229SH29, CUP: J53D23015950001). 

}
\section*{Conflict of interest/Competing interests}

{\small The author declares that he has no conflict of interest and competing interest.}

\section*{Availability of data and material and Code availability}

{ \small Not applicable.}

%

%%%%%%
\end{document}